\documentclass{amsart}
\usepackage{hyperref}
\usepackage[utf8]{inputenc}
\usepackage[cyr]{aeguill}
\usepackage[all]{xy}
\usepackage[french, english]{babel}
\usepackage[T1]{fontenc}
\usepackage{amscd}
\usepackage{amsmath}
\usepackage{amsthm} 
\usepackage{amssymb}
\usepackage{cleveref}
\usepackage{indentfirst}
\usepackage{enumitem}
\usepackage{array}
\newcolumntype{M}[1]{>{\center}p{#1}}
\newtheorem{theo}{Theorem}[section]
\newtheorem{theointro}{Theorem}

\newtheorem*{theo*}{Theorem}
\theoremstyle{definition}  
\newtheorem{defi}[theo]{Definition}
\newtheorem*{defi*}{Definition} 
                  
\newtheorem{lemm}[theo]{Lemma}
\newtheorem{prop}[theo]{Proposition}
\newtheorem{coro}[theo]{Corollary}

\newtheorem*{prop*}{Proposition}
\theoremstyle{remark}

\newtheorem{ex}[theo]{Example}
\newtheorem*{ex*}{Example}
\newtheorem{rema}[theo]{Remark}
\newtheorem*{rema*}{Remark}
\newtheorem*{Note*}{Notation}
\newtheorem{note}[theo]{Notation}
\newtheorem*{note*}{Notation}
\def\_{\underline}
\def\sh{\sharp}

\def\t{\mathfrak t}
\def\ı{\mathcal N}
\def\Op{\operatorname{Op}}

\def\Ind{\operatorname{Ind}}

\def\Q{\mathcal Q}
\def\Com{\operatorname{Com}}
\def\uCom{\operatorname{uCom}}
\def\MagCom{\operatorname{MagCom}}
\def\Lev{\operatorname{Lev}}
\def\D{\operatorname{D}}
\def\P{\mathcal{P}}
\def\M{\mathcal{M}}
\def\s{\mathfrak s}
\def\A{\mathcal{A}}
\def\Dp{\operatorname{D^{\pm}}}

\def\alg{\operatorname{alg}}

\def\val#1{\left\vert#1\right\vert}
\def\Set{\operatorname{Set}}
\def\Hom{\operatorname{Hom}}

\def\Coker{\operatorname{Coker}}

\def\Im{\operatorname{Im}}

\def\N{\mathbb N}
\def\diag{\shorthandoff{;:!?}
\xymatrix}

\def\0{\mathbf 0}
\def\Z{\mathbb Z}

\def\K{\scr K}

\def\U{\scr U}
\def\R{\mathbb R}
\def\F{\mathbb F}

\def\∑{\mathfrak S}
\def\to{\rightarrow}
\def\inj{\hookrightarrow}

\def\scr{\mathcal}

\def\val#1{{\vert#1\vert}}

\def\vect{\operatorname{vect}}

\def\mod{\operatorname{mod}}

\def\dar[#1]{\ar@<2pt>[#1]\ar@<-2pt>[#1]}
\usepackage[left=2cm, textwidth=16cm, textheight=22cm]{geometry}
\title{Unstable algebras over an operad}
\date{\today}
\begin{document}
\maketitle
\begin{center}
	\scshape Sacha Ikonicoff\footnote{Université de Paris, IMJ-PRG, CNRS, SU, 75013 Paris, France}
\end{center}
\begin{abstract}
  The aim of this article is to define and study a notion of unstable algebra over an operad that generalises the classical notion of unstable algebra over the Steenrod algebra. For this study we focus on the case of characteristic 2. We define $\star$-unstable $\P$-algebras, where $\P$ is an operad and $\star$ is a commutative binary operation in $\P$. We then build a functor that takes an unstable module $M$ to the free $\star$-unstable $\P$-algebra generated by $M$. Under some hypotheses on $\star$ and on $M$, we identify this unstable algebra as a free $\P$-algebra. Finally, we give some examples of this result, and we show how to use our main theorem to obtain a new construction of the unstable modules studied by Carlsson, Brown-Gitler, and Campbell-Selick, that takes into account their internal product.
\end{abstract}
\tableofcontents
\section{Introduction}
In this paper we study unstable modules over the mod 2 Steenrod algebra, taking as a starting point algebraic operations that appear naturally on certain unstable modules of interest.

The mod 2 Steenrod algebra $\A$ is the graded associative (non-commutative) unital $\F_2$-algebra generated by the degree $i$ elements $Sq^i$ for $i>0$, and satisfying relations called the Adem relations. An unstable $\A$-module is an $\A$-module satisfying the instability relation $Sq^ix=0$ for all $i>\val x$. An unstable algebra is a commutative, associative, graded algebra $(A,\cdot)$ endowed with a structure of unstable $\A$-module and satisfying two relations, one called the Cartan formula, and the other called the instability relation: $Sq_0 x=x\cdot x$, where $Sq_0x=Sq^{|x|}x$.

We use the formalism of (symmetric) algebraic operads to define and study algebraic operations on unstable modules. Following the classical definitions, we define a notion of unstable algebra over an operad $\P$ endowed with a commutative operation $\star$. Given such an operad $\P$ and a commutative operation $\star\in\P(2)^{\∑_2}$, a $\P$-algebra $M$ in unstable modules is said to be a $\star$-unstable $\P$-algebra if, for all $x\in M$, we have $Sq_0x=\star(x,x)$. We construct a functor $K_\P^{\star}$ that takes an unstable module $M$ to the free $\star$-unstable $\P$-algebra generated by $M$.

In a $\star$-unstable $\P$-algebra, the operation $\star$ satisfies an interchange law with respect to every operation of $\P$: for $\mu\in\P(n)$, $x_1,\dots,x_n\in M$, one has (see \Cref{*compalg}):
$$\star(\mu(x_1,\dots,x_n),\mu(x_1,\dots,x_n))=Sq_0\mu(x_1,\dots,x_n)=\mu(\star(x_1,x_1),\dots,\star(x_n,x_n)).$$
An operation $\star\in\P$ satisfying such a compatibility relation is called $\P$-central.

Our main result is the following Theorem:
\begin{theointro}[\Cref{theored}]
  Let $\P$ be an operad, $\star\in\P(2)^{\∑_2}$ be a $\P$-central operation, and $M$ be a connected reduced unstable module. There exists a graded $\P$-algebra isomorphism between the free $\star$-unstable $\P$-algebra $K_\P^\star(M)$ generated by $M$, and the free $\P$-algebra generated by $\Sigma \Omega M$.
\end{theointro}

The classical case of unstable algebras corresponds to the case of the $\cdot\,$-unstable $\Com$-algebras, where $\Com$ is the operad of (unital) commutative associative algebras, and $\cdot\in\Com(2)$ is the operadic generator. In this setting, the preceding result is known, notably for free unstable modules $M$ (see for example \cite{BK}). Moreover, when $M$ is freely generated by one element, the result of this theorem corresponds to the calculation due to Serre of the mod 2 cohomology of the Eilenberg-MacLane spaces of $\Z/2\Z$ \cite{JPS}.

The motivating examples for studying unstable algebras come from topology. Indeed, the mod 2 cohomology of a topological space is an unstable algebra. Other unstable modules are classically described with an internal product, such as Brown-Gitler modules, Carlsson modules, and Campbell-Selick modules. These unstable modules appear notably when studying the injective objects in the category of unstable modules. In \cite{DD}, Davis shows that the Carlsson module of weight 1, with the multiplication studied by Carlsson, which is commutative but not associative and satisfies an interchange law, is in fact isomorphic to the free `depth-invariant' algebra generated by an element of degree 1. His definition of a `depth-invariant' algebra corresponds to the definition of level algebras of Chataur and Livernet \cite{CL}.

A direct application of \Cref{theored} gives the following result:

\begin{prop*}[\Cref{propK1}]
	The Carlsson module $K(1)$ with its internal product is isomorphic to the free unstable level algebra generated by $F(1)$.
\end{prop*}

This assertion gives more precision to the result of Davis, which did not take into account the action of the Steenrod algebra.

So far, our result has not been generalised to the case of odd characteristic $p$. The definition of $\P$-central operations can be derived to operations of any arity $p$, and it is expected that the results extend to odd characteristic when considering unstable modules concentrated in even degree.

\begin{note}
  \item\begin{itemize}
    \item The base field is $\F:=\F_2$ for the whole article.
    \item $\F_{\vect}$ is the category of $\F$-vector spaces.
    \item $\A$ is the mod 2 Steenrod algebra. $\A_{\mod}$ is the category of left $\A$-modules.
    \item $\U$ is the category of unstable modules over the Steenrod algebra.
  \end{itemize}
\end{note}

\subsection*{Acknowledgements}

I want to thank Lionel Schwartz, whose ideas are the foundation of this work. I want to thank Geoffrey Powell, Christine Vespa and Muriel Livernet for fruitful conversations and valuable insight on the subject. Finally, I would like to thank the anonymous referee who suggested tremendous improvement in the layout of this article.

\section{Algebras over an operad in \texorpdfstring{$\U$}{U}}\label{algopU}
In this section, we recall the symmetric monoidal structure on the category of unstable modules over the Steenrod algebra, as well as the notion of symmetric operads. We focus on operads that come from operads in $\F_{\vect}$, concentrated in degree 0. We explain why algebras in the category of unstable modules over such an operad $\P$ are graded $\P$-algebras satisfying a generalised Cartan formula.

\begin{note}
  The mod 2 Steenrod algebra $\A$ is the associative, non-commutative, graded $\F_2$-algebra generated by the degree $i$ element $Sq^i$ for all $i>0$, satisfying the Adem relations:
  $$Sq^iSq^j-\sum_{k=0}^{\lfloor i/2\rfloor}\binom{j-k-1}{i-2k}Sq^{i+j-k}Sq^k=0,$$
  for all $i,j>0$ such that $i<2j$, where $\lfloor - \rfloor$ is the floor function, and where we denote by $Sq^0$ the unit of $\A$. We refer to \cite{LS} for notation and classical results about the Steenrod algebra.
  Recall that the category $\U$ is the full subcategory of $\A$-modules satisfying the condition: $Sq^j x=0$ for all $j>\val{x}$. The `top' square $x\mapsto Sq^{\val x}x$ is denoted by $Sq_0$.

  The category of $\A$-modules is endowed with a symmetric monoidal tensor product: if $M,N$ are two $\A$-modules, one can endow the graded tensor product $M\otimes N$ of the graded $\F$-vector spaces $M$ and $N$ with the $\A$-module structure given by:
  $$Sq^i(x\otimes y):=\sum_{k+l=i}Sq^kx\otimes Sq^ly,$$
  for all $x\in M$, $y\in N$, $i\in\N$. This monoidal structure actually comes from a Hopf algebra structure on $\A$ with cocommutative coproduct (see \cite{LS}, \cite{HM}).

  If $M$ and $N$ are unstable modules, then $M\otimes N$ is still unstable.
\end{note}
\subsection*{Recollections about operads}
We assume that the reader has a basic knowledge of operad theory in the algebraic setting. Our reference on the subject is the book \cite{LV} of Loday and Vallette. Let us recall the basic definition, in order to fix our notation:
\begin{defi}\item
  \begin{itemize}
    \item A symmetric sequence $\scr M$ is a sequence of vector spaces $(\scr M(n))_{n\in\N}$ such that, for all $n\in\N$, $\∑_n$ acts on $\scr M(n)$ on the right. The integer $n$ is often called ``arity''.
    \item Symmetric sequences form a category $\∑_{mod}$. This category is endowed with a tensor product such that:
    $$\left(\scr M\otimes\scr N\right)(n)=\bigoplus_{i+j=n}\Ind_{\∑_i\times\∑_j}^{\∑_n}\scr M(i)\otimes \scr N(j),$$
    where $\Ind_{\∑_i\times\∑_j}^{\∑_n}$ denotes the induced representation from the Young subgroup $\∑_i\times\∑_j$ of the group $\∑_n$.
    \item The category of symmetric sequences is endowed with another monoidal product $\circ$ given by:
    $$\left(\scr M\circ\scr N\right)(n)=\bigoplus_{k\ge0}\scr M(k)\otimes_{\∑_k} (\scr N^{\otimes k}(n)),$$
    with unit $\F$ concentrated in arity 1. For $\nu\in\M(n)$, and $\xi_1,\dots,\xi_n\in\scr N$, we denote by $(\nu;\xi_1,\dots,\xi_n)$ the element $\left[\nu\otimes \xi_1\otimes\dots\otimes \xi_n\right]_{\∑_n}\in\M\circ\scr N$.
    \item Operads are unital monoids in the monoidal category of symmetric sequences $(\∑_{mod},\circ)$. For an operad $\P$, we denote by $\gamma_{\P}:\P\circ\P\to\P$ its composition morphism, and $1_{\P}\in\P(1)$ its unit element. For $\nu\in\P(n)$, $\xi_1,\dots,\xi_n\in\P$, we denote by $\nu(\xi_1,\dots,\xi_n)$ the element $\gamma_\P(\nu;\xi_1,\dots,\xi_n)\in \P$. The partial compositions in an operad $\P$ are defined by 
    $$\mu\circ_i\nu=\mu(1_\P,\dots,1_\P,\underbrace{\nu}_{i\mbox{\scriptsize{-th input}}},1_\P,\dots,1_\P),$$ 
    for all $i\in\{1,\dots,m\}$, where $m$ is the arity of $\mu$.
    \item For any operad $\P$, $\P_{\alg}$ is the category of $\P$-algebras. We denote by $S(\P,-)$ the Schur functor $\F_{\vect}\to\F_{\vect}$ associated to $\P$, and the `free $\P$-algebra' functor $\F_{\vect}\to\P_{\alg}$, depending on the context. We recall that in $\F_{\vect}$, this functor is a monad and is defined by:
    $$S(\P,V)=\bigoplus_{n\ge0}\P(n)\otimes_{\∑_n}V^{\otimes n}.$$
    It can also be defined as $S(\P,V)=\P\circ V$ where $V$ is considered as a symmetric sequence concentrated in arity $0$.
    \item For an operad $\P$, a $\P$-algebra is an algebra over the monad $S(\P,-)$. In other terms, it is a couple $(V,\theta)$ where $V$ is a vector space and $\theta:S(\P,V)\to V$ is compatible with the composition and unit of $\P$. For $(V,\theta)$ a $\P$-algebra, $\mu\in\P(n)$ and $v_1,\dots,v_n\in V$, we denote by $\mu(v_1,\dots,v_n)$ the element $\theta(\mu;v_1,\dots,v_n)\in V$.
    \item $\Com$ (resp. $\uCom$) is the operad of commutative, associative algebras (resp. of commutative, associative, unital algebras).
  \end{itemize}
\end{defi}
\begin{defi}[see \cite{CL}]\label{defLEV}
A level algebra is a vector space $V$ endowed with a commutative bilinear operation $\star$ satisfying, for all $a,b,c,d\in V$, 
    $$(a\star b)\star(c\star d)=(a\star c)\star(b\star d).$$
	
  The operad of level algebras is denoted by $\Lev$. It is generated by an element $\star\in\Lev(2)^{\∑_2}$ satisfying the relation:
    $$\star(\star,\star)=\star(\star,\star) (2\ 3),$$
    where $(2\ 3)\in\∑_4$ is the transposition of $2$ and $3$. 
\end{defi}
\subsection*{Operads acting on unstable modules} The category of $\F$-vector spaces is identified with the subcategory of unstable modules concentrated in degree 0. This identification is compatible with the tensor products. For any operad $\P$ in $\F_{\vect}$, we deduce a notion of $\P$-algebra in $\A$, and in $\U$, regarding $\P$ as an operad concentrated in degree 0.
\begin{note}
  We denote by $\P_{\alg}^{\U}$ the category of $\P$-algebras in unstable modules. A morphism between $\P$-algebras in $\U$ is a $\P$-algebra morphism that is compatible with the action of $\A$.
\end{note}
\begin{prop}\label{propUCartan}
  A $\P$-algebra in $\U$ is a graded $\P$-algebra $M$ endowed with an action of the Steenrod algebra that satisfies the (generalised) Cartan formula, that is, for all $\mu\in\P(n)$, $(x_i)_{1\le i\le n}\in M^{\times n}$,
  $$Sq^i \mu(x_1,\dots,x_n)=\sum_{i_1+\dots+i_n=i}\mu(Sq^{i_1}x_1,\dots,Sq^{i_n}x_n).$$

  The forgetful functor $\P_{\alg}^{\U}\to \U$ admits, as a left adjoint functor, the functor $S(\P,-):\U\to\P_{\alg}^{\U}$, where, if $M$ is an unstable $\A$-module, $S(\P,M)$ is an unstable $\A$-module for the action induced by $M$ and the Cartan formula.
\end{prop}
\begin{proof}
  Since the elements of $\P$ are in degree $0$, for all $i>0$, $Sq^i$ acts trivially on the operadic elements. A $\P$-algebra in $\U$ is a graded $\P$-algebra $M$ endowed with an action of the Steenrod algebra such that the structural morphism $S(\P,M)\to M$ is compatible with the action of $\A$. The compatibility condition corresponds exactly to the Cartan formula (we will give an example of this computation in the case of the operad $\Com$). This also shows that $S(P,M)$, endowed with the action of the Steenrod algebra given by the Cartan formula, is the free $\P$-algebra in $\U$ generated by $M$.
\end{proof}
\begin{ex}
  All $\Com$-algebras in $\U$ satisfy the Cartan formula:

  Given an unstable module $M$, and $\theta:S(\Com,M)\to M$ a morphism in $\U$ endowing $M$ with the structure of a $\Com$-algebra in $\U$. Then $\theta$ is compatible with the action of the Steenrod algebra. Denote by $\cdot\in\Com(2)$ the operadic generator of $\Com$. Recall that the (associative and commutative) multiplication of $M$ is then defined by $(x,y)\mapsto\theta(\cdot;x,y)=\cdot(x,y)$. For all $x,y\in M$, $i\in\N$, one has:
   \begin{align*}
    Sq^i(\cdot(x,y))&=\sum_{j+k+l=i}(Sq^j\cdot)(Sq^kx,Sq^ly)\\
    &=\sum_{k+l=i}\cdot(Sq^kx, Sq^ly). 
  \end{align*}
\end{ex}
\begin{coro}\label{Sq0SP}
  The map $Sq_0$ is compatible with the action of $\P$. More precisely:
  
  Let $M$ be an unstable module, $\P$ an operad in $\F_{\vect}$. For all $\mu\in\P(k)$, $x_1,\dots,x_k\in M$, one has the following equality in $S(\P,M)$:
  \begin{equation*}
    Sq_0(\mu;x_1,\dots,x_k)=(\mu;Sq_0x_1,\dots,Sq_0x_k).
  \end{equation*}
\end{coro}
\begin{proof}
  Let $M$ be an unstable module and $\P$ be an operad in $\F_{\vect}$. For all $\mu\in \P(k)$, $x_1\in M^{n_1}$,\dots,$x_k\in M^{n_k}$, \Cref{propUCartan} gives:
  \begin{equation*}
    Sq^{n_1+\dots+n_k}(\mu;x_1,\dots,x_k)=\sum_{i_1+\dots+i_k=n_1+\dots+n_k}(\mu;Sq^{i_1}x_1,\dots,Sq^{i_k}x_k).
  \end{equation*}
  Because $M$ is an unstable module, one has $Sq^ix_j=0$ as soon as $i>n_j$. So the only terms of this sum that are not zero are the ones with $i_j\le n_j$ for all $j\in\{1,\dots,k\}$. The condition $i_1+\dots+i_k=n_1+\dots+n_k$ then implies that:
  \begin{equation*}
    Sq^{n_1+\dots+n_k}(\mu;x_1,\dots,x_k)=(\mu;Sq^{n_1}x_1,\dots,Sq^{n_k}x_k),
  \end{equation*}
  as required.
\end{proof}
\section{\texorpdfstring{$\P$}{P}-ideals in \texorpdfstring{$\U$}{U}}\label{ideauxU}
In this section, we recall the definition of a $\P$-ideal, where $\P$ is an operad, as well as the definition of the $\P$-ideal generated by a vector subspace of a $\P$-algebra. We then extend these definitions to those of a $\P$-ideal in $\U$ and of the $\P$-ideal in $\U$ generated by a vector subspace of a $\P$-algebra in $\U$. These objects have the desired universal properties in the corresponding category. 

\begin{defi}[\cite{GK}]
    \item\begin{itemize}
      \item Let $A$ be a $\P$-algebra. A $\P$-ideal of $A$ is a vector subspace $I$ of $A$ such that, for all $\mu\in\P(n)$, $a_1,\dots,a_n\in A$, 
      $$a_n\in I\ \Rightarrow\ \mu(a_1,\dots,a_n)\in I.$$
      The structure of $\P$-algebra of $A$ induces a structure of $\P$-algebra on the vector space $A/I$.
      \item Let $X\subset A$ be a vector subspace of the $\P$-algebra $A$. The $\P$-ideal generated by $X$, denoted by $(X)_\P$, is the smallest $\P$-ideal of $A$ that contains $X$. It satisfies the following universal property:
      For all $\P$-algebras $B$ and all $\P$-algebra morphisms $\varphi:A\to B$, $\varphi(X)=0$ if and only if $\varphi$ factors in a unique way into a $\P$-algebra morphism $\tilde \varphi:A/(X)_\P\to B$.
      \item A $\P$-ideal in $\U$ is a $\P$-ideal that is stable under the action of $\A$.
      \item Let $X\subset M$ be a vector subspace of the $\P$-algebra $M$ in $\U$. The $\P$-ideal in $\U$ generated by $X$, denoted by $(X)_{\P,\U}$, is the smallest $\P$-ideal in $\U$ of $M$ that contains $X$.
    \end{itemize}
\end{defi}
\begin{prop}\label{propideauxA}
  Let $M$ be a $\P$-algebra in $\U$, $N\subset M$ a sub-$\A$-module of $M$. Then $(N)_\P$ is stable under the action of $\A$. In particular, $(N)_{\P,\U}=(N)_\P$.
\end{prop}
\begin{proof}
  Every element of $(N)_\P$ is a sum of monomials of the type $\t=\mu(a_1,\dots,a_n)$, with $\mu\in \P(n)$, $a_1,\dots,a_{n-1}\in M$, and $a_n\in N$. For all $i\in\N$, one has:
  \begin{equation*}
      Sq^i\t=\sum_{i_1+\dots+i_n=i}\mu( Sq^{i_1}a_1,\dots, Sq^{i_{n-1}}a_{n-1},Sq^{i_n} a_n).
    \end{equation*}
    Yet, for all $i_n\in \N$, $Sq^{i_n}a_n \in N$, so $Sq^i\t\in(N)_{\P}$.

    In particular, $(N)_\P$ is a $\P$-ideal of $M$ in $\U$. For all $\P$-ideals $J$ of $M$ in $\U$ containing $N$, since $J$ is a $\P$-ideal containing $N$, one has $(N)_\P\subset J$. This proves that $(N)_\P=(N)_{\P,\U}$.
\end{proof}

\section{\texorpdfstring{$\star$}{star}-unstable \texorpdfstring{$\P$}{P}-algebras over the Steenrod algebra}\label{secalginst}
In this section, we define a notion of unstable algebra over the Steenrod algebra with respect to the data of an operad endowed with a commutative operation. This definition generalises the classical notion of unstable algebra over the Steenrod algebra. With this aim in mind, we recall the definition of the endofunctor $\Phi$ of the category of unstable modules, and we study the natural transformations linking $\Phi$ to $S(\P,-)$.

\begin{defi}[see \cite{LZ}, \cite{LS}]
  Denote by $\Phi:\U\to\U$ the functor defined on objects by:
  $$(\Phi(M))^n=\left\{\begin{array}{ll}
  M^{\frac n 2},&\mbox{ if }n\equiv 0\ [2],\\
  0,&\mbox{ otherwise.}
  \end{array}\right.$$
  For all $x\in M^n$, $\Phi x$ denotes the corresponding element in $(\Phi M)^{2n}$. For all $i\in\N$, one has:
  $$Sq^i \Phi x=\left\{\begin{array}{ll}
  \Phi(Sq^{\frac i 2} x),&\mbox{ if }i\equiv 0\ [2],\\
  0,&\mbox{ otherwise}.
  \end{array}\right.$$
  There is a natural transformation $\lambda:\Phi\to id_\U$ such that, for all $x\in M$, $\lambda_M(\Phi x)=Sq_0x$.
\end{defi}
\begin{prop}\label{lemimpairrepart}
	Let $\P$ be an operad, $M$ an unstable module. There is a natural inclusion $\Delta:\Phi\inj S^2$, where $S^2(M)=(M^{\otimes 2})_{\∑_2}$ is the second symmetric power, such that for all $x\in M$, $\Delta(\Phi x)=[x\otimes x]_{\∑_2}$. Moreover, for all $\star\in\P(2)^{\∑_2}$, there is a natural transformation $a^\star:S^2\to S(\P,-)$ that maps $[x\otimes y]_{\∑_2}\in S^2(M)$ to $(\star;x,y)\in S(\P,M)$. 
\end{prop}
\begin{proof}
	The only thing that is not clear is the compatibility of $\Delta$ with the action of the Steenrod algebra. 

	For all $x\in M$, one checks that
  \begin{align*}
    Sq^i\Delta(\Phi x)&=Sq^i[x\otimes x]_{\∑_2}=\sum_{j+k=i}[Sq^jx\otimes Sq^k x]_{\∑_2},\\
    &=\bigg(\sum_{j+k=i, j<k}[Sq^jx\otimes Sq^kx]_{\∑_2}+[Sq^kx\otimes Sq^j x]_{\∑_2}\bigg)+Y,
  \end{align*}
  where $Y:=\left\{\begin{array}{ll}
  [Sq^{\frac i 2}x\otimes Sq^{\frac i 2}x]_{\∑_2},&\mbox{ if }i\equiv 0\ [2],\\
  0,&\mbox{ otherwise.}
\end{array}\right.$. Yet, since $[Sq^jx\otimes Sq^kx]_{\∑_2}=[Sq^kx\otimes Sq^j x]_{\∑_2}$, one has $Sq^i\Delta(\Phi x)=[Sq^{\frac i 2}x\otimes Sq^{\frac i 2}x]_{\∑_2}=\Delta(\Phi(Sq^{\frac i 2}x) )=\Delta(Sq^i\Phi x)$.
\end{proof}
\begin{defi}\label{propalpha}
  Let $\P$ be an operad, and $\star\in\P(2)^{\∑_2}$. Denote by $\alpha^\star:\Phi\to S(\P,-)$ the composite $a^\star\circ \Delta$.

   A $\star$-unstable $\P$-algebra over the Steenrod algebra is a $\P$-algebra $(M,\theta)$ in $\U$ such that $\theta\circ\alpha^{\star}_M=\lambda_M$, that is, such that for all $x\in M$, 
 \begin{equation*}
 	Sq_0x=\star(x,x).
 \end{equation*}
  We denote by $\K_{\P}^{\star}$ the full subcategory of $\P_{\alg}^{\U}$ formed by $\star$-unstable $\P$-algebras. 
\end{defi}
\begin{ex}\item
	\begin{itemize}
		\item A $\cdot\,$-unstable $\uCom$-algebra is an unstable algebra in the classical sense (see \cite{SE}), that is, a $\uCom$-algebra $(M,\cdot)$ in $\U$ such that $Sq_0 x=x\cdot x$ for all $x\in M$.
		\item A $\star$-unstable $\Lev$-algebra is an unstable level algebra as defined in \cite{CL}, that is, a $\Lev$-algebra $(M,\star)$ in $\U$ such that $Sq_0 x=x\star x$ for all $x\in M$. Examples of such objects are given in \Cref{exapptheored}.
	\end{itemize}
\end{ex}
\begin{defi}\label{defislice}
  Denote by $\MagCom$ the operad freely generated by $\star\in \MagCom(2)^{\∑_2}$ and $(\MagCom\downarrow\Op)$ the category of operads under $\MagCom$ (see \cite{SML}).
\end{defi}
\begin{rema}
  An operadic morphism $\MagCom\to \P$ corresponds to the choice of a commutative operation $\star_\P\in\P(2)^{\∑_2}$. By abuse of notation we write $\star_\P$ for the corresponding object of $(\MagCom\downarrow\Op)$. With this notation, a morphism $f:\star_\P\to\star_\Q$ in $(\MagCom\downarrow\Op)$ corresponds to an operad morphism $f:\P\to\Q$ together with a choice $\star_\P\in\P(2)^{\∑_2}$, $\star_{\Q}\in\Q(2)^{\∑_2}$ such that $f(\star_\P)=\star_\Q$.

  Notice that any object $\star_\P$ in $(\MagCom\downarrow\Op)$ gives rise to the category $\K_\P^{\star}$ of $\star_\P$-unstable $\P$-algebras defined in \Cref{propalpha}.
\end{rema}
\begin{prop}\label{proprestrfunc}
  Let $\star_\P$, $\star_\Q$ be two objects of $(\MagCom\downarrow\Op)$. Any morphism $f:\star_\P\to\star_\Q$ naturally induces a restriction functor $f^*:\K_{\Q}^{\star_{\Q}}\to\K_{\P}^{\star_{\P}}$ which is the identity on the underlying unstable module.
\end{prop}
\begin{proof}
  An operad morphism $f:\P\to\Q$ induces a restriction functor $f^*:\Q_{\alg}\to\P_{\alg}$ taking a $\Q$-algebra $A$ to the $\P$-algebra with the same underlying vector space (see \cite{LV}) by setting, for all $\mu\in\P(n)$, $a_1,\dots,a_n\in A$,
  $$\mu(a_1,\dots,a_n) = f(\mu)(a_1,\dots,a_n).$$
  Let $\star_{\P}\in\P(2)^{\∑_2}$ and $\star_{\Q}=f(\star_{\P})$. If $A$ is a $\star_{\Q}$-unstable $\Q$-algebra, then $f^* A$ is clearly $\star_{\P}$-unstable.
\end{proof}
\section{\texorpdfstring{$\P$}{P}-central operations}\label{levcompsec}
In this section, we define the condition of $\P$-centrality for an operation $\star\in\P(2)^{\∑_2}$. Such an operation is said to be $\P$-central if it satisfies the interchange relation with respect to all other operations in $\P$. The condition of $\P$-centrality for an operation $\star\in\P(2)^{\∑_2}$ is necessary for the proof of \Cref{theored} in order to identify certain free $\star$-unstable $\P$-algebras.

\begin{note}
  In this section, the operad $\P$ in $\F_{\vect}$ is fixed.
\end{note}
\begin{defi}\label{deflevcomp}
  \item\begin{itemize}
    \item Let $\star \in\P(2)^{\∑_2}$. The operation $\star$ is said to be $\P$-central if it satisfies, for all $\mu\in\P(n)$, the interchange law:
    \begin{equation}\tag{E}
      \label{levcomp}\star (\mu,\mu)=\left(\mu(\underbrace{\star ,\dots,\star }_n)\right) \sigma_{2n},
    \end{equation}
    where $\sigma_{2n}\in\∑_{2n}$ maps $2i$ to $n+i$ and $2i-1$ to $i$ for all $i\in[n]$.
    \item Let $\star \in\P(2)^{\∑_2}$. A $\P$-algebra $A$ is said to be $\star$-compatible if, for all $\mu\in \P(n)$, $a_1,\dots,a_n\in A$, one has:
    $$\star(\mu(a_1,\dots,a_n),\mu(a_1,\dots,a_n))=\mu(\star(a_1, a_1),\dots,\star(a_n,a_n)).$$
  \end{itemize}
\end{defi}
\begin{rema}\label{levcomplev}
  A $\P$-central operation is a level operation. Indeed, relation (\ref{levcomp}) with $\mu=\star$ yields $\star(\star,\star)=\star(\star,\star) \sigma_{4}$. Since $\sigma_4=(2\ 3)$ is the transposition of $2$ and $3$, $\star$ is a level operation.
\end{rema}
\begin{rema}
  A $\P$-algebra $A$ is said to be $\star$-compatible if $\star$ acts on it as a $\P$-central operation, whether or not this operation satisfies the $\P$-centrality property. For instance, if $\star$ is $\P$-central, then all $\P$-algebras are $\star$-compatible, but an algebra on an operad $\P$ can be $\star$-compatible for an operation $\star\in\P(2)^{\∑_2}$ that is not $\P$-central.
\end{rema}
\begin{ex}
  Consider the operad $\MagCom$, which is freely generated by one commutative operation $\star\in\MagCom(2)^{\∑_2}$. Since $\star$ is not required to be a level operation, the equation (\ref{levcomp}) with $\star=\mu$ is not satisfied.

  It is clear that the multiplication of polynomials endows $\F[x]$ with a structure of $\MagCom$-algebra. It is easy to check that, with this definition, $\F[x]$ is $\star$-compatible.
\end{ex}
\begin{prop}\label{*compalg}
  All $\star$-unstable $\P$-algebras are $\star$-compatible.
\end{prop}
\begin{proof}
  Let $M$ be a $\star$-unstable $\P$-algebra, $\mu\in\P(n)$, and $x_1,\dots,x_n\in M$. One has
  \begin{align*}
    \star(\mu(x_1,\dots,x_n),\mu(x_1,\dots,x_n))&=Sq_0\mu(x_1,\dots,x_n),\\
    &=\mu(Sq_0x_1,\dots,Sq_0x_n),\\
    &=\mu(\star(x_1, x_1),\dots,\star(x_n,x_n)),
  \end{align*}
  where the second equality is a consequence of \Cref{Sq0SP}.
\end{proof}
The next proposition allows one to check, from a presentation of an operad $\P$, if an operation $\star\in\P(2)^{\∑_2}$ is $\P$-central. It will be used in \Cref{Levcompop}.
\begin{prop}\label{Levcompgen}
   Let $\star\in\P(2)^{\∑_2}$. Let $F$ be a sub-symmetric sequence of $\P$. Suppose that $F$ generates the operad $\P$.

  Then $\star$ is $\P$-central if and only if it satisfies relation (\ref{levcomp}) of \Cref{deflevcomp} for all $\mu\in F$.
\end{prop}
\begin{proof}
  It suffices to check that if $\mu,\nu\in F$ satisfy (\ref{levcomp}), then all partial compositions $\mu\circ_i\nu$ satisfy (\ref{levcomp}), and, in this setting, if $\mu$ and $\nu$ have same arity, $\mu+\nu$ satisfies (\ref{levcomp}). The additivity statement is easily checked. Let $\mu,\nu\in F$. Denote by $m$ and $n$ the respective arities of $\mu$ and of $\nu$. If $\mu$ and $\nu$ satisfy (\ref{levcomp}), then:
  \begin{align*}
    \star(\mu\circ_i\nu,\mu\circ_i\nu)&=\left(\left(\star(\mu,\mu)\right)\circ_{m+i}\nu\right)\circ_i\nu,\\
    &=\left(\left(\mu(\star,\dots,\star) \sigma_{2m}\right)\circ_{m+i}\nu\right)\circ_i\nu,\\
    &=\left(\left(\mu(\star,\dots,\star)\circ_{i+1}\nu\right)\circ_{i}\nu\right) \sigma',\\
    &=\left(\mu\left(\star,\dots,\star,\underbrace{\star(\nu,\nu)}_i,\star,\dots,\star\right)\right) \sigma',\\
    &=\left(\mu\left(\star,\dots,\star,\underbrace{\nu(\star,\dots,\star) \sigma_{2n}}_i,\star,\dots,\star\right)\right) \sigma',\\
    &=(\mu\circ_i\nu)(\star,\dots,\star) \sigma_{2(m+n-1)},
  \end{align*}
  where $\sigma'\in\∑_{2m+2n-2}$ is the block permutation obtained by applying $\sigma$ to $2m$ blocks of size $1$, except for the $i$-th and the $(i+1)$-th, of size $n$.
\end{proof}
Let us now give some examples of operads endowed with central operations.
\begin{ex}
	The generator of the operad $\Com$ is $\Com$-central. The generator of the operad $\Lev$ defined in \ref{defLEV} is $\Lev$-central. 
\end{ex}
\begin{defi}
  \item\begin{itemize}
    \item Set $\D=\F[d]$, the polynomial algebra in one indeterminate $d$, seen as an operad concentrated in arity $1$, with unit $1\in\F$. A $\D$-algebra is a vector space endowed with an endomorphism $d$.
    \item For all $s>0$, set $Q_s \D=\F[d]/\left(d^s-1\right)$, seen as an operad concentrated in arity $1$. A $Q_s\D$-algebra is a vector space endowed with an endomorphism $d$ such that $d^s$ is the identity morphism.
    \item For all $q>0$, set $T_q \D=\F[d]/\left(d^{q+1}\right)$, seen as an operad concentrated in arity $1$. A $T_q\D$-algebra is a vector space endowed with a nilpotent endomorphism $d$ of order $\le q+1$.
    \item Set $\Dp=\F[d,d^{-1}]$, the Laurent polynomial in one indeterminate $d$, seen as an operad concentrated in arity $1$. A $\Dp$-algebra is a vector space endowed with an automorphism $d$.
    \item Denote by $T_q\Lev$ the $q$-th truncation of the operad $\Lev$, that is, the quotient of $\Lev$ by the operadic ideal spanned by all compositions of the generator such that the composition tree (see \cite{LV}) is of height greater than $q$.
  \end{itemize}
\end{defi}
\begin{lemm}
  Let $\P$ be an operad. The following morphism of $\∑$-modules is a distributive law (see \cite{LV}):
  \begin{eqnarray*}
  &\D\circ\P&\to \P\circ\D\\
     &(d^n;\mu)&\mapsto (\mu; \underbrace{d^n,\dots,d^n}_{m})
  \end{eqnarray*}
  and it induces a distributive law on $\P\circ Q_sD$ and on $\P\circ \Dp$.
\end{lemm}
\begin{proof}
  It is a straightforward verification.
\end{proof}
\begin{lemm}\label{Levcompop}
  The operadic generator $\star\in\Lev(2)^{\∑_2}$ yields a $T_q\Lev$-central operation. The operation $(\cdot;d,d)$ is a $\uCom\circ\D$-central operation, a $\uCom\circ\Dp$-central operation, and a $\uCom\circ Q_s\D$-central operation. More generally, if $\star\in\P(2)$ is $\P$-central, then $(\star;d^i,d^i)$ is $\P\circ\D$, $\P\circ\Dp$, and $\P\circ Q_s\D$-central, for all $i\in\N$.
\end{lemm}
\begin{proof}
  All these assertions are proved by use of \Cref{Levcompgen}.
\end{proof}
The following statement will be useful in \Cref{exapptheored}.
\begin{prop}\label{opmorph}
  The operadic morphism $\Lev\to \Com\circ\D$ (respectively, $\Lev\to \Com\circ\Dp$), mapping the generator $\star\in\Lev(2)^{\∑_2}$ to $(\cdot;d,d)$ is a monomorphism. It passes to the quotients, inducing a monomorphism $T_q\Lev\to \Com\circ T_q\D$.
\end{prop}
\begin{proof}
	The vector spaces $\Lev(n)$ are endowed with a basis constituted of ordered partitions $J=(J_0,\dots,J_p)$ of the set $\{1,\dots,n\}$, $p$ being any integer, satisfying $\sum_{i=0}^p\frac{|J_i|}{2^i}$, and where we allow some of the $J_i$'s to be empty (see \cite{SI}).

	The vector spaces $\Com\circ\D(n)$ are endowed (for $n>0$) with the same basis without the summation condition $\sum_{i=0}^p\frac{|J_i|}{2^i}$.

	The operadic morphism $\Lev\to \Com\circ\D$, mapping the generator $\star\in\Lev(2)^{\∑_2}$ to $(\cdot;d,d)$ injects the given basis of $\Lev$ into the given basis of $\Com\circ\D$.

	The vector space $T_q\Lev(n)$ (respectively, $\Com\circ T_q\D(n)$) is the quotient vector space of $\Lev(n)$ (resp. $\Com\circ\D(n)$) by the vector subspace spanned by the partitions $J=(J_0,\dots,J_p)$ such that there exist an $i>q$ satisfying $J_i\neq \emptyset$. Hence, the inclusion of bases passes to the quotients.
\end{proof}
\section{\texorpdfstring{$\star$}{*}-unstable \texorpdfstring{$\P$}{P}-algebras generated by an unstable module}\label{secalginsteng}
In this section, we build, for any operad $\P$ endowed with a commutative operation $\star\in\P(2)^{\∑_2}$, a functor that assigns to an unstable module $M$ the free $\star$-unstable $\P$-algebra generated by $M$.

In the beginning of this section, we recall some basic notions for the study of unstable modules. We refer to \cite{LS} for this notation. 

The main result of this section is \Cref{theored}, which gives a concise description of the free $\star$-unstable $\P$-algebra generated by a connected reduced unstable module, when the operation $\star$ satisfies the $\P$-centrality condition defined in Section \ref{levcompsec}. Under these conditions, we identify this $\star$-unstable $\P$-algebra, which is a quotient of a $\P$-algebra, as a free $\P$-algebra.

When working with the operad $\uCom$, endowed with its multiplication $\cdot\in\uCom(2)^{\∑_2}$, and when we consider free unstable algebras generated by a free, monogeneous unstable module, the result of \Cref{theored} corresponds to the calculation due to Serre of the mod 2 cohomology of the Eilenberg-MacLane spaces of $\Z/2\Z$ \cite{JPS}.

\begin{note}
  Throughout this section, we fix an operad $\P$ and a commutative operation $\star\in\P(2)^{\∑_2}$.
\end{note}
\begin{defi}[see \cite{SE}, \cite{LS}, \cite{LZ}]\label{defmiscunst}
  \item\begin{itemize}
    \item The suspension functor $\Sigma:\U\to\U$ takes an unstable module $M$ to the unstable module $\Sigma M$ defined by $(\Sigma M)^d=M^{d-1}$, and $Sq^i(\sigma x)=\sigma(Sq^ix)$, where $\sigma x\in M^{d+1}$ corresponds to the element $x\in M^d$.
   
    \item The functor $\Sigma$ admits a left adjoint denoted by $\Omega$. For all $M\in\U$, the unstable module $M/\Im(Sq_0)$ is a suspension, and the de-suspension $\Sigma^{-1}(M/\Im(Sq_0))$ is isomorphic to $\Omega M$. The counit of the adjunction $(\Omega,\Sigma)$ is a natural isomorphism $\Omega\Sigma M\cong M$.

    \item Let $I=(i_1,\dots,i_k)$ be a (finite) sequence of integers. Then $I$ is called admissible if for all $h\in\{1,2,\dots,k-1\}$, $i_h\ge 2i_{h+1}$. The excess of an admissible sequence is the integer $e(I):=i_1-i_2-\dots-i_k$.

    \item Let $I=(i_1,\dots,i_k)$ be a (finite) sequence of integers. We denote by $Sq^I$ the product $Sq^{i_1}\dots Sq^{i_k}$ in $\A$.
  \end{itemize}
\end{defi}
\begin{defi}[see \cite{LS}]
  For $n\in\N$, the free unstable module generated by one element $\iota_n$ in degree $n$ is denoted $F(n)$. One has $\Hom_\U(F(n),M) \cong M^n$.

  The following assertions are all consequences of the definition of $F(n)$:
  \begin{itemize}
  \item \label{lemFn} The unstable module $F(n)$ is isomorphic to:
    $$\Sigma^n\A/(Sq^I:I \mbox{ admissible, and }e(I)>n).$$
    Therefore, $F(n)$ admits as a graded vector space basis the set of $Sq^I\iota_n$ where $\iota_n$ is the generator of degree $n$ and $I$ satisfies $e(I)\le n$. In particular, the unstable module $F(1)$ has for basis the set $\{j_k\}_{k\in\N}$, where $j_k:=Sq^{2^{k-1}}\dots Sq^1\iota_1\in F(1)^{2^k}$ ($j_0=\iota_1$).
    \item\label{omegaFn} For all $n>0$, $\Omega F(n)$ is isomorphic as an unstable module to $F(n-1)$. Indeed, for $M$ an unstable module, there is a one-to-one correpondence (natural in $M$):
    $$\Hom_\U(\Omega F(n),M)\cong \Hom_\U(F(n),\Sigma M)\cong \left(\Sigma M\right)^n\cong M^{n-1}\cong\Hom_\U(F(n-1),M).$$
\end{itemize}
\end{defi}
\begin{defi}[see \cite{LS}]\item
 \begin{itemize}
   \item An unstable module $M$ is said to be reduced if the application $Sq_0:M\to M$ is injective.
  
    \item An unstable module $M$ is said to be connected if $M^0=0$.
 \end{itemize}
 \end{defi}
\begin{note}
  \item\begin{itemize}
    \item For all $k\in\N$, the element $\star_k\in\P(2^k)$ is inductively defined by:
  $$\star_0=1_\P,\mbox{ and, }\forall k\ge 1,\ \star_k=\star(\star_{k-1},\star_{k-1}).$$
  \item For all $\mu\in\P(n)$, $x\in M$, where $M$ is an unstable module, the element $(\mu;\underbrace{x,\dots,x}_{n})\in S(\P,M)$ is denoted by$(\mu;x^{\times n})$. If $M$ is a $\P$-algebra, the element $\mu(\underbrace{x,\dots,x}_{n})\in M$ is denoted by $\mu(x^{\times n})$.
  \end{itemize}
\end{note}
\begin{lemm}
  If $\star$ is $\P$-central, then $\star_k$ belongs to $\P(2^k)^{\∑_{2^k}}$. 
\end{lemm}
\begin{proof}
  Let us show this result by induction. For $k=0$, this is obvious since $\∑_1$ is trivial.

  Suppose that, for a given $k\in\N$, $\star_k$ is stable under the action of $\∑_k$. The compatibility between operadic composition and symmetric group action implies that $\star_{k+1}$ is stable under the action of the wreath product $\∑_2\wr\∑_{2^k}\subset \∑_{2^{k+1}}$. Furthermore, assuming that $\star$ is $\P$-central, one checks easily that the element $\sigma_{2^{k+1}}\in\∑_{2^{k+1}}$ defined in \ref{deflevcomp} stabilises $\star_{k+1}$. 

  It suffices to show that $\∑_2\wr\∑_{2^k}\cup\{\sigma_{2^{k+1}}\}$ generates $\∑_{2^{k+1}}$. One can show that every transposition $(i\ j)$ is in the subgroup $G < \∑_{2^{k+1}}$ generated by $\∑_2\wr\∑_{2^k}\cup\{\sigma_{2^{k+1}}\}$. Note that if $i,j\in\{1,\dots,2^k\}$ or $i,j\in\{2^k+1,\dots,2^{k+1}\}$, then $(i\ j)\in \∑_{2^k}\times\∑_{2^k}\subset \∑_2\wr\∑_{2^k}$. Let us suppose that $i\in\{1,\dots,2^k\}$ and $j\in\{2^k+1,\dots,2^{k+1}\}$. Note that $\sigma_{2^{k+1}}(1)=1$ and $\sigma_{2^{k+1}}(2^{k+1})=2$. Then $\sigma_{2^{k+1}}^{-1}(1\ 2)\sigma_{2^{k+1}}$ is the transposition $(1\ 2^k+1)$. So, $(1\ i)(2^k+1\ j)\sigma_{2^{k+1}}^{-1}(1\ 2)\sigma_{2^{k+1}}(1\ i)(2^k+1\ j)=(i\ j)$. This shows that $(i\ j)\in G$.
\end{proof}
\begin{prop}
  Let $\P$ be an operad, $\star\in\P(2)^{\∑_2}$ be a commutative operation and $M$ be an unstable module. We set:
  $$K_\P^\star(M):=S(\P,M)/(\{Sq_0\t+\star(\t,\t):\t\in S(\P,M)\})_{\P,\U}.$$
  Then $K_\P^\star(M)$ is a $\star$-unstable $\P$-algebra. Moreover, $K_\P^\star:\U\to\K_\P^\star$ gives a left adjoint functor for the forgetful functor $U:K_\P^\star\to\U$.
\end{prop}
\begin{proof}
  Let $\t\in S(\P,M)$. One has $Sq_0[\t]=[Sq_0\t]=[\star(\t,\t)]=\star([\t],[\t])$ in $K_\P^\star(M)$, so $K_\P^\star(M)$ is $\star$-unstable.

  Let $N$ be a $\star$-unstable $\P$-algebra, and $g:M\to N$ a morphism of $\A$-modules. Since $S(\P,M)$ is the free $\P$-algebra in $\U$ generated by $M$, there exists a unique morphism $g':S(\P,M)\to N$ of $\P$-algebras in $\U$ that extends $g$.

  Now, as $K_\P^\star(M)$ is $\star$-unstable, $g'((\{Sq_0\t+\star(\t,\t):\t\in S(\P,M)\})_{\P,\U})=0$. So there is a unique factorisation morphism $g'':K_\P^\star(M)\to N$. Thus, since $K_\P^\star(M)$ is $\star$-unstable, it is the free $\star$-unstable $\P$-algebra generated by $M$.
\end{proof}
The following proposition, using the notation introduced in \Cref{secalginst}, will be useful in \Cref{exapptheored}.
\begin{prop}\label{propmorphrestr}
  Let $\star_\P$, $\star_\Q$ be two objects in $(\MagCom\downarrow\Op)$, and $M$ be an unstable module. Any morphism $f:\star_\P\to\star_\Q$ naturally induces a morphism $f_*:K_\P^{\star_\P}(M)\to f^* K_\Q^{\star_\Q}(M)$, where $f^*:\K_\Q^{\star_\Q}\to \K_\P^{\star_\P}$ is the restriction functor from \Cref{proprestrfunc}.
\end{prop}
\begin{proof}
  The morphism $S(f,M):S(\P,M)\to f^* S(\Q,M)$ passes to the quotients, inducing the desired morphism $f_*:K_\P^{\star_\P}(M)\to f^* K_\Q^{\star_\Q}(M)$.
\end{proof}
\begin{lemm}[see \cite{LS}]
  Let $M$ be an unstable module. The unstable modules $\Sigma\Omega M$ and $\Coker\lambda_M$ are isomorphic. The following diagram is a short exact sequence when $M$ is reduced:
  $$\diag{\Phi M\ar[r]^-{\lambda_M}&\ar[r]^-{pr}M&\Sigma\Omega M}.$$
  Moreover, the morphism $pr:M\to \Sigma\Omega M$ is the unit of the adjunction $(\Omega,\Sigma)$.
\end{lemm}
\begin{defi}\label{secglob}
Let $M$ be a reduced unstable module. A graded section, denoted by $s:\Sigma\Omega M\to M$, is the data, for all $d\in\N$, of a linear section $s:(\Sigma\Omega M)^{d}\to M^{d}$ of the map $pr:M^d\to(\Sigma\Omega M)^d$ (that is, such that $pr\circ s=id_{\Sigma\Omega M}$). We draw the reader’s attention to the fact that a graded section is not, in general, compatible with the action of $\A$, but is only a graded linear map.
\end{defi}
\begin{theo}\label{theored}
  Let $\P$ be an operad in $\F_{\vect}$, $\star\in\P(2)^{\∑_2}$ be a $\P$-central operation. For all connected reduced unstable modules $M$, there exists an isomorphism of graded $\P$-algebras between the $\star$-unstable $\P$-algebra $K_\P^\star(M)$ and the free $\P$-algebra generated by $\Sigma\Omega M$.

  This isomorphism is not natural in $M$, and depends on the choice of a graded section $s:\Sigma\Omega M\to M$ (see \Cref{secglob}).
\end{theo}
\begin{proof}
  We refer to Section \ref{preuvered}.
\end{proof}

\section{Proof of \texorpdfstring{\Cref{theored}}{}}\label{preuvered}
Throughout this section, we fix an operad $\P$ in $\F_{\vect}$, endowed with a $\P$-central operation $\star\in\P(2)^{\∑_2}$, and a connected unstable module $M$ endowed with a graded section $s:\Sigma\Omega M\to M$ (see \Cref{secglob}). For the proof of \Cref{preuvered3} and \Cref{preuvered1}, we do not assume that $M$ is reduced.

The $\P$-algebra $K_\P^\star(M)$ is described as a quotient of the free $\P$-algebra $S(\P,M)$ by an ideal. In \Cref{preuvered3}, \Cref{preuvered1} and \Cref{preuvered2}, we simplify this ideal, replacing it by a $\P$-ideal generated by a vector subspace built from $\Sigma\Omega M$.

We then give a construction for the desired $\P$-algebra isomorphism between $K_\P^{\star}(M)$ and the free $\P$-algebra $S(\P,\Sigma\Omega M)$.

\begin{note}\label{notepreuvered}
  Let $Unst=\{Sq_0x+(\star;x,x):x\in M\}\subset S(\P,M)$, and $X=\{Sq_0\t+(\star;\t,\t):\t\in S(\P,M)\}\subset S(\P,M)$.
\end{note}
\begin{lemm}\label{preuvered3}
  The vector subspace $X\subset S(\P,M)$ defined in \ref{notepreuvered} is stable under the action of $\A$. In particular, one deduces from \Cref{propideauxA} that:
  $$K_\P^{\star}(M)=S(\P,M)/(X)_{\P}.$$
\end{lemm}
\begin{proof}
 Using the notation of \Cref{secalginst}, there are natural transformations $\lambda:\Phi\to id_{\U}$ and $\alpha^\star:\Phi\to S(\P,-)$. These yield natural transformations $\lambda_{S(\P,-)}=\Phi\circ S(P,-)\to S(P,-)$ and $\gamma_\P\circ\alpha^\star_{S(\P,-)}:\Phi\circ S(\P,-)\to S(\P,-)$. The vector subspace $X$ is the image of the map $\lambda_{S(\P,M)}+\gamma_\P\circ\alpha^\star_{S(\P,M)}$, which is compatible with the action of $\A$. This implies that $X$ is stable under the action of $\A$.
\end{proof}
\begin{prop}\label{preuvered1}
In $S(\P,M)$, one has $(X)_{\P}=(Unst)_{\P}$, where $X$ and $Unst$ are defined in \ref{notepreuvered}. In particular, one deduces from \Cref{preuvered3} that:
$$K_\P^{\star}(M)=S(\P,M)/(Unst)_{\P}.$$
\end{prop}
\begin{proof}
  Let us show that $Unst\subset X$ and $X\subset (Unst)_{\P}$. The first inclusion is clear, let us prove the second one. Let $\t:=(\mu;x_1,\dots,x_n)\in S(\P,M)$ be a $\P$ monomial, with $x_1,\dots,x_n\in M$. Following \Cref{Sq0SP}, one has:
    $$Sq_0\t=(\mu;Sq_0x_1,\dots,Sq_0x_n),$$
  So the following element is in $(Unst)_{\P}$:
  \begin{align*}
    Sq_0\t+\mu((\star;x_1,x_1),\dots,(\star;x_n,x_n))&=Sq_0\t+(\mu(\star,\dots,\star);x_1,x_1,\dots,x_n,x_n),\\
    &=Sq_0\t+(\mu(\star,\dots,\star)  \sigma_{2n};x_1,\dots,x_n,x_1,\dots,x_n).
  \end{align*}
  Since $\star$ is $\P$-central, this element is equal to:
  $$Sq_0\t+(\star(\mu,\mu);x_1,\dots,x_n,x_1,\dots,x_n)=Sq_0\t+\star(\t,\t).$$
  Hence, $X\subset (Unst)_{\P}$.
\end{proof}
\begin{rema}
  Note that the $\P$-centrality hypothesis is crucial. The mildest hypothesis needed on $\star\in\P(2)^{\∑_2}$ for the proof of \Cref{preuvered1} is that for all $\mu\in\P(n)$, if $V$ denotes the vector space spanned by $\{x_1,\dots,x_n\}$, then in $S(\P,V)$, one has
    $$(\mu(\star,\dots,\star)  \sigma_{2n};x_1,\dots,x_n,x_1,\dots,x_n)=(\star(\mu,\mu);x_1,\dots,x_n,x_1,\dots,x_n),$$
    and this is equivalent to the $\P$-centrality hypothesis.
\end{rema}
\begin{lemm}\label{preuveredpre4}
  Assume that $M$ is reduced. Then, for all $n\in \N$,
  $$M^n=\{Sq_0^ks(b):b\in(\Sigma\Omega M)^{2^{-k}n},k\in\N\},$$
  where we set $(\Sigma\Omega M)^{2^{-k}n}=0$ when $2^{-k}n\notin\N$.
\end{lemm}
\begin{proof}
  This is a reformulation of the fact that $M$ being reduced, it is freely generated by $s(\Sigma\Omega M)$ under the action of $Sq_0$.
\end{proof}
\begin{defi}\label{defisplitred}
  Set:
  $$E:=\left\{Sq_0^ks(b)+\left(\star_{k};\left(s(b)\right)^{\times 2^{k}}\right): b\in \Sigma\Omega M, k\in\N\right\}\subset S(\P,M).$$
\end{defi}
\begin{lemm}\label{preuvered2}
  Suppose that $M$ is reduced. One has $(E)_\P=(Unst)_\P$, where $E$ is defined in \ref{defisplitred} and $Unst$ is defined in \ref{notepreuvered}. One then deduces from \Cref{preuvered1} that:
  $$K_\P^\star(M)=S(\P,M)/(E)_{\P}.$$ 
\end{lemm}
\begin{proof}
  Let us show that $E\subset (Unst)_\P$ and $Unst\subset (E)_\P$.

  \vspace{.5cm}
  \textbf{Proof of the inclusion $Unst\subset (E)_\P$.}

  Let $x\in M$. According to \Cref{preuveredpre4}, $x$ can be written as a sum of elements of the form $y:=Sq_0^{k}s(b)\in M$ where $k\in\N$, $b\in\Sigma\Omega M$. 

  On the one hand, $Sq_0y=Sq_0^{k+1}s(b)$, so the following element $\alpha$ is in $(E)_\P$:
  $$\alpha:=Sq_0y+\left(\star_{k+1};(s(b))^{\times 2^{k+1}}\right)$$

  On the other hand, $(\star;y,y)=\left(\star;Sq_0^{k}s(b),Sq_0^{k}s(b)\right)$. So the following element is in $(E)_\P$:
  $$\beta:=(\star;y,y)+\star\left(\left(\star_{k};(s(b))^{\times 2^{k}}\right),\left(\star_{k};(s(b))^{\times 2^{k}}\right)\right)=(\star;y,y)+\left(\star_{k+1};(s(b))^{\times 2^{k+1}}\right).$$

  Thus, the elements of the form $\alpha+\beta=Sq_0y+(\star;y,y)$
  are in $(E)_\P$, so $Unst\subset(E)_\P$.

  \vspace{.5cm}
  \textbf{Proof of the inclusion $E\subset (Unst)_\P$.}

  Let $k\in\N$, $b\in \Sigma\Omega M$. Let us show, by induction on $k\in\N$, that $Sq_0^ks(b)+\left(\star_k;\left(s(b)\right)^{\times2^k}\right)\in(Unst)_\P$. For $k=0$, one has $Sq_0^ks(b)+\left(\star_k;\left(s(b)\right)^{\times2^k}\right)=0\in(Unst)_\P$.

  Let us assume that, for all $l\in\N$, one has $Sq_0^ls(b)+\left(\star_l;\left(s(b)\right)^{\times2^l}\right)\in(Unst)_\P$.
  Let $k=l+1$. Note that the following element is in $Unst$:
  $$\alpha:=Sq_0^ks(b)+(\star;Sq_0^ls(b),Sq_0^ls(b)).$$
  By the induction hypothesis, the following element is in $(Unst)_\P$:
  $$Sq_0^ls(b)+\left(\star_l;(s(b))^{\times2^{l}}\right).$$
  Thus, the following element is in $(Unst)_\P$:
  $$\beta:=\star\left(Sq_0^ls(b)+\left(\star_l;(s(b))^{\times2^{l}}\right),Sq_0^ls(b)+\left(\star_l;(s(b))^{\times 2^{l}}\right)\right)=\left(\star;Sq_0^ls(b),Sq_0^ls(b)\right)+\left(\star_k;(s(b))^{\times 2^{k}}\right).$$
  Finally, the following element is in $(Unst)_\P$:
  \begin{equation*}
    \alpha+\beta=Sq_0^ks(b)+\left(\star_k;\left(s(b)\right)^{\times2^k}\right).\qedhere
  \end{equation*}
\end{proof}
\begin{theo}[\Cref{theored}]
  Suppose that $M$ is reduced. the graded section $s$ induces a $\P$-algebra isomorphism $K_\P^\star(M)\cong S(\P,\Sigma\Omega M)$.
\end{theo}
\begin{proof}
  Recall from \Cref{preuvered2} that $K_\P^\star(M)\cong S(\P,M)/(E)_\P$, where $E$ is defined in \ref{defisplitred}.

  The graded section $s:\Sigma\Omega M\to M$ induces a graded $\P$-algebra monomorphism $S(\P,s):S(\P,\Sigma\Omega M)\to S(\P,M)$, which, in turn, induces a graded $\P$-algebra map $\psi_s:S(\P,\Sigma\Omega M)\to S(\P,M)/(E)_\P$.

  Conversely, using \Cref{preuveredpre4}, one can define a graded linear map $\varphi_s:M\to S(\P,\Sigma\Omega M)$ by setting:
  $$\varphi_s(Sq_0^kb)=\star_k(b^{\times 2^k}),$$
  for all $b\in\Sigma\Omega M$, $k\in\N$. This map induces a graded $\P$-algebra map $\bar\varphi_s:S(\P,M)\to S(\P,\Sigma\Omega M)$.

  Let us show that $\bar \varphi_{s}$ factorises into a $\P$-algebra morphism $\hat \varphi_{s}:S(\P,M)/(E)_\P\to S(\P,\Sigma\Omega M)$. Since $\bar\varphi_{s}$ is compatible with the action of $\P$, and since $\bar\varphi_{s}(s(b))=b$ for all $b\in\Sigma\Omega M$, one has
  $$\bar\varphi_{s}(Sq_0^ks(b))=\left(\star_k;b^{\times 2^k}\right)=\left(\star_k\left(\left(\overline\varphi_{s}(s(b))\right)^{\times 2^k}\right)\right)=\overline\varphi_{s}\left(\star_k;\left(s(b)\right)^{\times 2^k}\right).$$

  We have proved the existence of the factorisation morphism $\hat \varphi_{s}:S(\P,M)/(E)_\P\to S(\P,\Sigma\Omega M)$. 

  Let us show that $\psi_s$ is a bijection, and that $\hat\varphi_s$ is its inverse.

  The vector space $S(\P,\Sigma\Omega M)$ is spanned by elements of the form $(\nu;b_1,\dots,b_m)$, where $\nu\in \P(n)$ and $b_1, \dots, b_m\in \Sigma\Omega M$. But, one has:
  \begin{equation*}
    \hat\varphi_s\circ\psi_s(\nu;b_1,\dots,b_m)=\hat \varphi_{s}(\nu;s(b_1),\dots,s(b_m))=(\nu;b_1,\dots,b_m).
  \end{equation*}
  So $\hat\varphi_{s}\circ\psi_s=id_{S(\P,\Sigma\Omega M)}$, and $\psi_s$ is injective.

  To complete the proof, it suffices to show that $\psi_s$ is surjective. For this purpose, let $\pi:S(\P,M)\to S(\P,M)/(E)_{\P}$ denote the canonical projection, so that $\pi\circ S(\P,s)=\psi_s$. Since $\pi$ is onto, it suffices to prove that for all $\s\in S(\P,M)$, there is an $\s'\in S(\P,\Sigma \Omega M)$ such that $\psi_s(\s')=\pi(\s)$ .

 \Cref{preuveredpre4} implies that $S(\P,M)$ is spanned by the element of the type:
  $$\s=(\nu;Sq_0^{k_{1}}s(b_{1}),\dots,Sq_0^{k_{m}}s(b_{m})),$$
 where $\nu\in \P(m)$, $k_1,\dots,k_m\in\N$, and $b_{1},\dots,b_{m}\in\Sigma\Omega M $. This element can also be written:
 \begin{multline*}
    \s=\underbrace{(\nu;Sq_0^{k_{1}}s(b_{1}),\dots,Sq_0^{k_{m}}s(b_{m}))+\left(\nu(\star_{k_{1}},\dots,\star_{k_{m}});\left(s(b_{1})\right)^{\times 2^{k_{1}}},\dots,\left(s(b_{m})\right)^{\times 2^{k_{m}}}\right)}_{\in (E)_\P}\\
  +\underbrace{\left(\nu(\star_{k_{1}},\dots,\star_{k_{m}});\left(s(b_{1})\right)^{\times 2^{k_{1}}},\dots,\left(s(b_{m})\right)^{\times 2^{k_{m}}}\right)}_{\in\Im\left(S(\P,s)\right)}.
  \end{multline*}

  Let $\s'$ be the following element of $S(\P,\Sigma\Omega M)$:
  $$\s'=\left(\nu(\star_{k_{1}},\dots,\star_{k_{m}});b_{i}^{\times 2^{k_{1}}},\dots,b_{m}^{\times 2^{k_{m}}}\right).$$
  We have found an $\s'$ such that $\pi(\s)=\psi_s(\s')$. So $\psi_s$ is onto, hence it is a bijection, and $\hat\varphi_s$ is its inverse.
  \end{proof}
  
\section{First examples and applications}\label{firstex}
In this section, we study some applications of \Cref{theored}, when we take $M$ to be the free unstable module $F(n)$. \Cref{theored} gives, under some assumptions on the unstable module $M$ and the operation $\star\in\P(2)^{\∑_2}$, an isomorphism of $\P$-algebras $K_\P^{\star}(M)\cong S(\P,\Sigma\Omega M)$. This isomorphism is not natural in $M$, and depends on a graded section $s$ of $M$ (see \Cref{secglob}). In the case where $M=F(1)$, we will see that there is a unique choice for $s$. When $M=F(n)$, we will give a somewhat natural choice. In these cases, we will study the action of $\A$ obtained on $S(\P,\Sigma\Omega M)$ by transfer from the action of $\A$ on $K_\P^\star(M)$. We will then show that the result of \Cref{theored} needs not hold when $M$ is not reduced.

We will use the notation used in the proof of \Cref{theored}. The isomorphism $K_\P^{\star}(M)\cong S(\P,\Sigma\Omega M)$ will be denoted by $\hat\varphi_s$. It is the factorisation of a graded $\P$-algebra morphism $\bar\varphi_s:S(\P,M)\to S(\P,\Sigma\Omega M)$.

\begin{rema}\label{remhatvarphi}
	Let $\P$ be an operad in $\F_{\vect}$, and let $M$ be a connected reduced unstable module endowed with a graded section $s$. Following the constructions of the proof of \Cref{theored}, the inverse of $\hat\varphi_{s}:K_\P^{\star}(M)\to S(\P,\Sigma\Omega M)$ is the $\P$-algebra morphism sending $x\in\Sigma\Omega M$ to $[s(x)]_{E}\in K_\P^\star(M)$.
\end{rema}
\begin{defi}
	Let $\P$ be an operad in $\F_{\vect}$, $M$ be a connected reduced unstable module, endowed with a graded section $s$. One defines an action of $\A$ on the graded $\P$-algebra $S(\P,\Sigma\Omega M)$ by setting:
	$$Sq^i\odot \t:=\hat\varphi^{-1}(Sq^i\hat\varphi(\t)).$$
\end{defi}
\begin{defi}
      \item\begin{itemize}
        \item The classical graded section $s:\Sigma F(n-1)\to F(n)$ sends $\sigma(Sq^I\iota_{n-1})\in \Sigma F(n-1)$ to $Sq^I\iota_n$.
        \item In the case $n=1$, since $\Sigma F(0)$ only contains one non-zero element $\sigma\iota_0$ of degree $1$. Since $F(1)^1$ only contains one non-zero element $\iota_1$, the only graded section $\Sigma F(0)\to F(1)$ sends $\sigma\iota_0$ to $\iota_1$.
       \end{itemize} 
  \end{defi}
\begin{prop}\label{propKF1}
  $K_\P^\star(F(1))$ is the free $\P$-algebra generated by one element $\iota_1$ of degree $1$ endowed with the unstable action of $\A$ determined by:
  $$Sq^j\iota_1:=\left\{\begin{array}{ll}
    \iota_1,&\mbox{if }j=0,\\
    \star(\iota_1,\iota_1),&\mbox{if }j=1,\\
    0,&\mbox{otherwise,}
  \end{array}\right.$$
 and satisfying the Cartan formula.
 \end{prop}
 \begin{proof}
   Let us describe the action of $\A$ obtained on the graded $\P$-algebra $S(\P,\Sigma F(0))$ from the action of $\A$ on $K_\P^\star(F(1))$ through the isomorphism $\hat\varphi_{s}$ deduced from $s$ as defined above. Since $\hat\varphi_{s}$ is an isomorphism of $\P$-algebras, and since $K_\P^\star(F(1))$ satisfies the Cartan formula, it suffices to describe the action of $\A$ on the generator $\iota_1$. Since $Sq^0\odot\iota_1=\iota_1$, and since $Sq^i\odot\iota_1=0$ for all $i\neq0,1$, it suffices to compute $Sq^1\iota_1$. But, because $K_\P\star(F(1))$ is $\star$-unstable, one necessarily gets $Sq^1\odot\iota_1=(\star;\iota_1,\iota_1)$.
 \end{proof}
\begin{ex}
    This example shows that for $x\in \Omega M$, even when $i<\val{x}$, the equality $Sq^i\odot(\sigma x)=\sigma(Sq^ix)$ does not hold in $S(\P,\Sigma\Omega M)$.

    Set $M=F(2)$, $x:=\sigma(Sq^4Sq^2Sq^1\iota_1)\in \Sigma F(1)\subset S(\P,\Sigma F(1))$. Adem relations give $Sq^1Sq^4=Sq^5$. So, the action of $\A$ on $\Sigma F(1)$ gives $Sq^1x=\sigma(Sq^5Sq^2Sq^1\iota_1)$, which is $0$ because of the instability condition. However, the action of $K_\P^{\star}(F(2))$ transfered to $S(\P,\Sigma F(1))$ through the isomorphism $\hat\varphi_{s}$, with $s$ as defined above yields:
    \begin{align*}
      Sq^1\odot x&=\bar\varphi(Sq^1s(x))=\bar\varphi(Sq^5Sq^2Sq^1\iota_2)=\bar\varphi(Sq_0Sq^2Sq^1\iota_2)=\bar\varphi(\star;Sq^2Sq^1\iota_2,Sq^2Sq^1\iota_2)\\&=\bar\varphi(\star;s(\sigma Sq^2Sq^1\iota_1),s(\sigma Sq^2Sq^1\iota_1))=(\star;\sigma Sq^2Sq^1\iota_1,\sigma Sq^2Sq^1\iota_1).
    \end{align*}
    If $\star\neq 0$ in $\P$, one deduces that $Sq^1\odot(\sigma x)\neq\sigma(Sq^ix)$ in $S(\P,\Sigma\Omega M)$.
\end{ex}
\begin{lemm}
	Let $\P$ be an operad in $\F_{\vect}$, $\star\in\P(2)^{\∑_2}$ be a $\P$-central operation. Let $M$ be a connected reduced unstable module. For all graded sections $s:\Sigma\Omega M\to M$, the action of $K_{\P}^{\star}(M)$ transferred to $S(\P,\Sigma\Omega M)$ through the $\P$-algebra isomorphism $\hat\varphi_{s}$ always yields:
  $$Sq_0\odot(\mu;x_1,\dots,x_n)=\star\left((\mu;x_1,\dots,x_n),(\mu;x_1,\dots,x_n)\right),$$
  where $\mu\in \P(n)$ in $x_1,\dots,x_{n}\in \Sigma \Omega M$.
\end{lemm}
\begin{proof}
  It is a consequence of the $\star$-instability of $K_{\P}^{\star}(M)$. More precisely, one has, with the notation of the proof of \Cref{theored},
  \begin{align*}
    Sq_0\odot(\mu;x_1,\dots,x_n)&=\bar\varphi(Sq_0\mu(s(x_1),\dots,s(x_n)))\\
    &=\bar\varphi(\star(\mu(s(x_1),\dots,s(x_n)),\mu(s(x_1),\dots,s(x_n)))).
  \end{align*}
  Since $\bar\varphi$ is a $\P$-algebra morphism, one deduces that:
  \begin{align*}
    Sq_0\odot(\mu;x_1,\dots,x_n)&=\star\left(\bar\varphi(\mu(s(x_1),\dots,s(x_n))),\bar\varphi(\mu(s(x_1),\dots,s(x_n)))\right).
  \end{align*}
  The result then follows from \Cref{remhatvarphi}
\end{proof}

Let us now show that, when the generating unstable module is not reduced, the conclusion of \Cref{theored} needs not hold.
\begin{ex}
  Set $M=\Sigma F(0)$, $\P=\uCom$ and $\star=\cdot$, the generator of $\uCom$. Recall that $\Sigma F(0)$ is a dimension $1$ vector space concentrated in degree $1$ with generator $\sigma \iota_0$ satisfying $Sq_0(\sigma\iota_0)=Sq^1\sigma\iota_0=0$. Hence, it is not reduced. In the unstable algebra $K_{\uCom}^\cdot(\Sigma F(0))$, one has $\cdot(\sigma\iota_0,\sigma\iota_0)=Sq^1\sigma\iota_0=0$. So $K_{\uCom}^\cdot(\Sigma F(0))$ is isomorphic to $\Sigma F(0)$ endowed with a trivial multiplication. On the other hand, recall (see \Cref{omegaFn}) that $\Omega\Sigma F(0)\cong F(0)$, so $S(\uCom,\Sigma\Omega (\Sigma F(0)))\cong S(\uCom,\Sigma F(0))$ is the polynomial algebra on one degree one element. It is clear that, as vector spaces, $K^\cdot_{\uCom}(\Sigma F(0))$ is not isomorphic to $S(\uCom,\Sigma\Omega\Sigma F(0))$.
\end{ex}
\section{Further applications}\label{exapptheored}
In this section, we recall the definition and structures of several classical unstable modules, such as Brown-Gitler modules, Carlsson modules, and Campbell-Selick modules. These modules come equipped with internal products satisfying different properties. We show that some of these classical modules, with their operations, are free unstable algebras over operads.

We refer to \cite{LZ} and \cite{LS} for the definitions and results on Brown-Gitler modules, the Brown-Gitler algebra, Carlsson modules and the Carlsson algebra, and we refer to \cite{CS} for the definitions and results on Campbell-Selick modules.

For a different approach to unstable modules, the reader may refer to \cite{K1} and \cite{K3}, where the author gives a presentation of the Brown-Gitler and Carlsson modules using the theory of generic representations.

\begin{note}
	For a vector space $V$, we denote by $V^{\sh}=\Hom_{\F}(V,\F)$ the dual vector space.
\end{note}

\begin{defi}[\cite{LZ}]
\item\begin{itemize}
	\item Let $n\in\N$. The $n$-th Brown-Gitler module $J(n)$ is a representing object of the functor $H^n:\U\to\Set$, mapping $M$ to $(M^n)^\sh$. For all $n,m>0$, there is a linear correspondence $J(n)^m\cong (F(m)^n)^\sh$.

	\item The unstable modules $J(n)$ are endowed with an external product $\mu_{m,n}:J(n)\otimes J(m)\to J(n+m)$. The map $\mu_{m,n}\in\Hom_{\U}(J(n)\otimes J(m),J(n+m))\cong\left(\left(J(n)\otimes J(m)\right)^{n+m}\right)^\sh$ is the only non-zero element.
	\item The unstable module $J(n)$ is endowed with an internal product obtained as a composite:
	$$\diag{J(n)\otimes J(n)\ar[r]^-{\mu_{n,n}}& J(2n)\ar[r]^{\cdot Sq^n}& J(n)},$$
	where the second map $\cdot Sq^n\in \Hom_{\U}(J(2n),J(n))\cong\left(J(2n)^{n}\right)^\sh\cong F(n)^{2n}$ corresponds to the element $Sq_0\iota_n\in F(n)^{2n}$.
	\item The direct sum $J:=\bigoplus_{n\in\N}J(n)$, with the multiplication given by the external products $\mu_{m,n}$, is the Brown-Gitler algebra, also called the Miller algebra. An element of $J(n)$, seen in $J$, is said to have weight $n$, and this weight is additive with respect to the multiplication of $J$.
\end{itemize}
\end{defi}
\begin{theo}[Miller \cite{HM}]\label{theoBG}
The Brown-Gitler algebra $J$ is isomorphic to the polynomial algebra $\F[x_i,i\in\N]$, with $\val{x_i}=1$, and $x_i$ has weight $2^i$, endowed with the unstable action of $\A$ induced by the `shifted action':
\begin{equation}\tag{SA}\label{shac}
	Sq^jx_{i}:=\left\{\begin{array}{ll}
    x_i,&\mbox{if }j=0,\\
    x_{i-1}^2,&\mbox{if }j=1,\\
    0,&\mbox{otherwise,}
  \end{array}\right.
\end{equation}
  where we set $x_{-1}=0$, and satisfying the Cartan formula.
\end{theo}
\begin{defi}[\cite{GC},\cite{LZ}]
	\item\begin{itemize}
		\item The $n$-th Carlsson module $K(n)$ is the limit of the following diagram in $\U$:
	$$\diag{J(n)&J(2n)\ar[l]_{\cdot Sq^n}&\cdots\ar[l]&\ar[l]J(2^qn)&J(2^{q+1}n)\ar[l]_{\cdot Sq^{2^q}n}&\cdots\ar[l]}.$$
	The external products $\mu_{n,m}$ and the internal products on $J(n)$ pass to the limits, yielding an external product $\mu_{n,m}:K(n)\otimes K(m)\to K(n+m)$ and an internal product $K(n)\otimes K(n)\to K(n)$.
		\item The direct sum $K:=\bigoplus_{n\in\N}K(n)$, with the multiplication given by the external products $\mu_{m,n}$, is the Carlsson algebra.
	\end{itemize}
\end{defi}
\begin{theo}[Carlsson \cite{GC}, \cite{LZ}]\label{theoCarl}
  The Carlsson algebra $K$ is isomorphic to the polynomial algebra $\F[x_i,i\in\Z]$, with $\val{x_i}=1$, endowed with the unstable action of $\A$ induced by the shifted action (\ref{shac}) and satisfying the Cartan formula.
\end{theo}
\begin{defi}
	The $s$-th Campbell-Selick module $M_s$ is the polynomial algebra $\F[x_i,i\in\Z/s\Z]$, with $\val{x_i}=1$, endowed with the unstable action of $\A$ induced by the shifted action (\ref{shac}) and satisfying the Cartan formula.
\end{defi}
\begin{rema}
  It is clear from the definition of the shifted action (\ref{shac}) that the Brown-Gitler algebra, the Carlsson algebra and the Campbell-Selick modules are not unstable algebras in the classical acceptation of the term. For instance, if they were unstable, one would have $Sq_0x_i=Sq^1x_i=x_i^2$, which is not the case.
\end{rema}
\subsection*{Applications of our results by identification of structures}
  We have seen that the Brown-Gitler algebra, the Carlsson algebra and the Campbell-Selick modules are not unstable algebras in the classical sense. However, one can check that the Brown-Gitler algebra $J$ is a $(\cdot;d,d)$-unstable $\uCom\circ\D$-algebra, that the Carlsson algebra $K$ is a $(\cdot;d,d)$-unstable $\uCom\circ\Dp$-algebra, that $J(n)$ and $K(n)$ with their internal products are unstable $\Lev$-algebras, and that $M_s$ is a $(\cdot;d,d)$-unstable $Q_s\D$-algebra. 

  We will now identify some of these modules as free unstable algebras over their respective operads.

  We refer to the notation introduced in \Cref{levcompsec} for operads and their central operations.

\begin{ex}
	The free $\cdot\,$-unstable $\uCom$-algebra generated by $F(1)$ is isomorphic to $H^*(\R P^\infty,\F)$ with the cup-product (see \cite{SE}).

	The free $\cdot\,$-unstable $\uCom$-algebra generated by $F(n)$ is isomorphic to the $\uCom$-algebra generated by $\Sigma F(n-1)$, and this correspond to the description of the mod 2 cohomology of the $n$-th Eilenberg-MacLane space of $\Z/2\Z$ given in \cite{JPS}.
\end{ex}
\begin{prop}\label{propK1}
	The Carlsson module $K(1)$ with its internal product is isomorphic to the free unstable level algebra generated by $F(1)$.
\end{prop}
\begin{proof}
	Following Davis's result (\cite{DD}, Theorem 4.2), $K(1)$ with its internal product is the free level algebra on one generator of degree 1. This corresponds to the underlying level algebra structure of $K_{\Lev}^{\star}(F(1))$. Since $K$ is endowed with the unstable action of $\A$ induced by the shifted action (\ref{shac}) and satisfying the Cartan formula, \Cref{propKF1} shows that $K(1)$ with its internal product is isomorphic to $K_{\Lev}^{\star}(F(1))$ as unstable level algebras.
\end{proof}
\begin{prop}\label{propK}\item
	\begin{itemize}
		\item The Carlsson algebra $K$, with the algebra endomorphism $d:K\to K$ mapping $x_i$ to $x_{i-1}$, is the free $(\cdot;d,d)$-unstable $\uCom\circ\Dp$-algebra generated by $F(1)$.
		\item The Campbell-Selick module $M_s$, with its algebra structure and the algebra endomorphism $d:M_s\to M_s$ mapping $x_i$ to $x_{i-1}$, is the free $(\cdot;d,d)$-unstable $\uCom\circ Q_s\D$-algebra generated by $F(1)$.
		\item The direct sum of Brown-Gitler modules $\bigoplus_{i=1}^{2^q}J(i)$ with the multiplication given by the external product $\mu_{m,n}$, and the endomorphism $d:\bigoplus_{i=1}^{2^q}J(i):\bigoplus_{i=1}^{2^q}J(i)$ defined on $J(i)$ by
    $$d=\begin{cases}
      \cdot Sq^{i/2}&\mbox{ if }i\equiv 0\ [2],\\
      0&\mbox{ otherwise,}
    \end{cases}$$
    is the free $(\cdot;d,d)$-unstable $\uCom\circ T_q\D$-algebra generated by $F(1)$.
	\end{itemize}
\end{prop}
\begin{proof}\item
	\begin{itemize}
		\item \Cref{theoCarl} show that the unitary commutative and associative algebra $K$ is free on $\{x_i\}_{i\in\Z}$, of degree $1$. With the (invertible) map $d$ that sends $x_i$ to $x_{i-1}$, $K$ is endowed with a structure of $\uCom\circ\Dp$-algebra structure, and it is clearly freely generated by $x_0$. The result is then deduced by comparing the action of $\A$ on $K$ given in \Cref{theoCarl} with \Cref{propKF1}.
		\item The proof is similar. As a unitary commutative and associative algebra, $M_s$ is defined as free on $\{x_i\}_{i\in\Z/s\Z}$. With the (cyclic) map $d$ that sends $x_i$ to $x_{i-1}$, $M_s$ is endowed with a structure of $\uCom\circ Q_s\D$-algebra structure, and it is clearly freely generated by $x_0$. The result is then deduced by comparing the action of $\A$ on $M_s$ with \Cref{propKF1}.
		\item The proof is again similar. As a unitary commutative and associative algebra, \Cref{theoBG} shows that $\bigoplus_{i=1}^{2^q}J(i)$ is, as a sub-algebra of $J$, freely generated by $\{x_i\}_{0\le i\le q}$. With the map $d$ (nilpotent of degree $q+1$), $\bigoplus_{i=1}^{2^q}J(i)$ is clearly freely generated by $x_q$. The result is then deduced by comparing the action of $\A$ on $J$ given in \Cref{theoBG} with \Cref{propKF1}.
	\end{itemize}
\end{proof}
\subsection*{Morphisms of operads, and links between associated unstable algebras}Recall from \Cref{proprestrfunc} and \Cref{propmorphrestr} that any  morphism $f:\star_\P\to \star_\Q$ in $(\MagCom\downarrow \Op)$ induces a restriction functor $f^*:\K_\Q^{\star_\Q}\to \K_\P^{\star_\P}$, and a morphism $f_*:K_\P^{\star_\P}(M)\to f^* K_\Q^{\star_\Q}(M)$ in $\K_\P^\star$, natural with respect to the unstable module $M$.
\begin{note}
  By abuse of notation, we will drop the restriction functor and denote by $K_\Q^{\star_\Q}(M)$ the $\star_\P$-unstable algebra $f^* K_\Q^{\star_\Q}(M)$.
\end{note}
\begin{ex}
  The monomorphism $\Lev\to\uCom\circ\Dp$ mapping $\star\in\Lev(2)^{\∑_2}$ to $(\cdot;d,d)\in (\uCom\circ\Dp)(2)^{\∑_2}$ induces a monomorphism of unstable level algebras $K_{\Lev}^{\star}(F(1))\to K_{\uCom\circ\Dp}^{(\cdot;d,d)}(F(1))$. This morphism corresponds to the inclusion $K(1)\inj K$ of the Carlsson module of weight 1 into the Carlsson algebra.
\end{ex}
\begin{prop}\label{propBGTqLev}
  The Brown-Gitler module of weight $2^q$, $J(2^q)$, with the multiplication given by the internal product $(\cdot Sq^{2^q})\circ\mu_{2^q,2^q}$ is the free $\star$-unstable $T_q \Lev$-algebra generated by $F(1)$.
\end{prop}
\begin{proof}
  The monomorphism $T_q\Lev\to \uCom\circ T_q\D$ induces a monomorphism of $\star$-unstable $T_q\Lev$-algebras $K_{T_q\Lev}^{\star}(F(1))\to K_{\uCom\circ T_q\D}^{(\cdot;d,d)}(F(1))$. This morphism corresponds to the inclusion of the free $T_q\Lev$-algebra generated by $x_q$ as a subalgebra of $\bigoplus_{i=1}^{2^q}J(i)$. It is easy to show that the image of this inclusion is $J(2^q)$ (for example, by adapting the proof of \cite{DD}, Theorem 4.2).
\end{proof}
\begin{prop}
  The truncation tower of $\Lev$ induces the following cofiltration of $K(1)$ in the category of unstable level algebras:
  $$\diag{J(1)&J(2)\ar[l]_{\cdot Sq^1}&\ar[l]\cdots&J(2^q)\ar[l]&J(2^{q+1}\ar[l]_{\cdot Sq^{2^q}})&\cdots\ar[l]}.$$
\end{prop}
\begin{proof}
  The operad $\Lev$ is the limit of the diagram:
  $$\diag{T_0\Lev&T_1\Lev\ar[l]&\cdots\ar[l]&T_q\Lev\ar[l]&T_{q+1}\Lev\ar[l]&\cdots\ar[l],}$$
  where the morphisms $T_{q+1}\Lev\to T_q\Lev$ are the operadic projections. Following \Cref{propmorphrestr}, the previous diagram induces a diagram:
  $$\diag{K_{T_0\Lev}^{\star}(F(1))&K_{T_1\Lev}^{\star}(F(1))\ar[l]&\cdots\ar[l]&K_{T_q\Lev}^{\star}(F(1))\ar[l]&K_{T_{q+1}\Lev}^{\star}(F(1))\ar[l]&\cdots\ar[l],}$$
  in the category of unstable level algebras. For all $q\in\N$, \Cref{propBGTqLev} gives the isomorphism $K_{T_q\Lev}^{\star}(F(1))\cong J(2^q)$. It is easy to check that the morphism $K_{T_{q+1}\Lev}^{\star}(F(1))\to K_{T_q\Lev}^{\star}(F(1))$ induced by the operadic projection corresponds to $\cdot Sq^{2^q}:J(2^{q+1})\to J(2^q)$. Indeed, $\cdot Sq^{2^q}$ is compatible with the internal products on $J(2^{q+1})$ and $J(2^q)$, and maps $x_{q+1}$ to $x_q$. The unstable module $K(1)$ is defined as the limit of the $J(2^q)$'s under the morphisms $\cdot Sq^{2^q}$, and its internal product is given as the limit of the internal products of the $J(2^q)$'s. Hence, $K(1)$ with its internal product is the limit of the desired diagram in the category of unstable level algebras.
\end{proof}

\end{document}